\documentclass[11pt]{amsart2000}
\newtheorem{thm}{Theorem}[section]
\newtheorem{cor}[thm]{Corollary}

\newtheorem{note}[thm]{Note}

\newtheorem{res}[thm]{Result}

\theoremstyle{definition}
\newtheorem{eg}[thm]{Example}
\newtheorem{defn}[thm]{Definition}

\newcommand{\piccap}[1]{\stepcounter{figure}
\smallskip
Figure \thefigure #1}

\newcommand{\Z}{\mathbb Z}
\newcommand{\N}{\mathbb N}
\newcommand{\Q}{\mathbb Q}
\newcommand{\R}{\mathbb R}

\begin{document}
\setlength{\unitlength}{0.01in}
\linethickness{0.01in}
\begin{center}
\begin{picture}(474,66)(0,0)
\multiput(0,66)(1,0){40}{\line(0,-1){24}}
\multiput(43,65)(1,-1){24}{\line(0,-1){40}}
\multiput(1,39)(1,-1){40}{\line(1,0){24}}
\multiput(70,2)(1,1){24}{\line(0,1){40}}
\multiput(72,0)(1,1){24}{\line(1,0){40}}
\multiput(97,66)(1,0){40}{\line(0,-1){40}}
\put(143,66){\makebox(0,0)[tl]{\footnotesize Proceedings of the Ninth Prague Topological Symposium}}
\put(143,50){\makebox(0,0)[tl]{\footnotesize Contributed papers from the symposium held in}}
\put(143,34){\makebox(0,0)[tl]{\footnotesize Prague, Czech Republic, August 19--25, 2001}}
\end{picture}
\end{center}
\vspace{0.25in}
\setcounter{page}{135}
\title{Sub-representation of posets}
\author{M. K. Gormley}
\thanks{The research of the first author was supported by a distinction
award scholarship from the Department of Education for Northern Ireland}
\author{T. B. M. McMaster}
\address{Department of Pure Mathematics\\
Queen's University, Belfast\\
University Road\\
Belfast, BT7 1NN\\
United Kingdom}
\email{m.k.gormley@qub.ac.uk}
\email{t.b.m.mcmaster@qub.ac.uk}
\thanks{M. K. Gormley and T. B. M. McMaster,
{\em Sub-representation of posets},
Proceedings of the Ninth Prague Topological Symposium, (Prague, 2001),
pp.~135--146, Topology Atlas, Toronto, 2002}
\begin{abstract}
We define a property {\it sub-representability} and we give a complete
characterisation of sub-representability of posets. 
\end{abstract}
\maketitle

\section{Introduction}

\begin{defn}
Given an ordered set $E$ and a topological space $X$, we say that $E$ can 
be {\it realised within $X$} [see \cite{realise1}] if there is an 
injection $j$ from $E$ into the class of (homeomorphism classes of) 
subspaces of $X$ such that, for $x,y$ in $E$, $x\leq y$ if and only if
$j(x)$ is homeomorphically embeddable into $j(y)$.
\end{defn}

The question of which spaces have the `converse' of this property of {\it
realisability} appears difficult in general, but we are able to handle the
principal $T_0$ case. What we shall now do therefore is to try to
represent the family of sub-posets of a partially ordered set $P$, ordered
by embeddability, within that family ordered by inclusion. However, we
additionally want the representation to be such that each sub-poset is
represented by an embeddability-wise equivalent poset. That is, we wish to
sub-represent the poset.

\begin{defn}
Whenever $P_1$ and $P_2$ are posets, we shall say that $P_1\hookrightarrow 
P_2$ if
and only if $P_1$ is isomorphic to a subset of $P_2$.
\end{defn}

\begin{defn}\label{subrep}
We shall say that $P$ is {\it sub-representable} if and only if there
exists a map $g$ from ${\mathbb P}(P)$ to ${\mathbb P}(P)$ such that
for all
$P_1, P_2
\in {\mathbb P}(P)$:
\begin{itemize}
\item[(i)]
$P_1\hookrightarrow P_2$ if and only if $g(P_1)\subseteq g(P_2)$, and
\item[(ii)]
$P_1\hookrightarrow g(P_1)\hookrightarrow P_1$.
\end{itemize}
\end{defn}

\begin{eg}
\addtocounter{figure}{1}
Figure \thefigure \;details a sub-representation of a four-point poset.
\addtocounter{figure}{-1}
\begin{center}
\begin{picture}(250,400)(0,0)
\put(50,40){\circle*{10}}
\put(0,80){\circle*{10}}
\put(0,120){\circle*{10}}
\put(10,110){\vector(1,2){45}}
\put(0,85){\line(0,1){30}}
\put(0,180){\circle*{10}}
\put(0,185){\line(0,1){30}}
\put(0,220){\circle*{10}} 
\put(0,225){\line(0,1){30}}
\put(0,260){\circle*{10}} 
\put(90,100){\circle*{10}}
\put(70,100){\circle*{10}} 
\put(80,200){\circle*{10}} 
\put(70,240){\circle*{10}}
\put(80,205){\line(-1,4){7.6}}
\put(80,205){\line(1,4){7.6}}
\put(90,240){\circle*{10}}
\put(50,300){\circle*{10}}
\put(50,340){\circle*{10}}
\put(30,380){\circle*{10}}
\put(70,380){\circle*{10}}
\put(50,305){\line(0,1){30}} 
\put(50,345){\line(1,2){16}}
\put(50,345){\line(-1,2){16}}
\put(40,50){\vector(-1,1){25}}
\put(60,50){\vector(2,3){20}}
\put(0,140){\vector(0,1){20}} 
\put(10,267){\vector(1,1){27}}
\put(80,250){\vector(-1,2){20}}
\put(80,120){\vector(0,1){60}}
\put(220,40){\circle*{10}} \put(200,35){$3$}
\put(170,80){\circle*{10}} \put(150,75){$2$} 
\put(170,120){\circle*{10}} \put(150,115){$3$}
\put(180,110){\vector(1,2){45}}
\put(170,85){\line(0,1){30}} \put(170,180){\circle*{10}}
\put(150,175){$1$} \put(170,185){\line(0,1){30}}
\put(170,220){\circle*{10}} \put(150,215){$2$}
\put(170,225){\line(0,1){30}} \put(170,260){\circle*{10}}
\put(150,255){$3$} \put(260,100){\circle*{10}} \put(275,95){$4$}
\put(240,100){\circle*{10}} \put(220,95){$3$}
\put(250,200){\circle*{10}} \put(265,195){$2$}
\put(240,240){\circle*{10}} \put(220,235){$3$}
\put(250,205){\line(-1,4){7.6}} \put(250,205){\line(1,4){7.6}}
\put(260,240){\circle*{10}} \put(275,235){$4$}
\put(220,300){\circle*{10}} \put(235,295){$1$}
\put(220,340){\circle*{10}} \put(235,335){$2$}
\put(200,380){\circle*{10}} \put(180,375){$3$}
\put(240,380){\circle*{10}} \put(255,375){$4$}
\put(220,305){\line(0,1){30}} \put(220,345){\line(1,2){16}}
\put(220,345){\line(-1,2){16}} \put(210,50){\vector(-1,1){25}}
\put(230,50){\vector(2,3){20}} \put(170,140){\vector(0,1){20}}
\put(180,267){\vector(1,1){27}} \put(250,250){\vector(-1,2){20}}
\put(250,120){\vector(0,1){60}}
\put(10,15){\scriptsize Sub-posets ordered}
\put(10,3){\scriptsize by embeddability.}
\put(170,15){\scriptsize Representative sub-posets}
\put(170,3){\scriptsize ordered by inclusion.}
\end{picture}\\
\piccap{} 
\end{center}
\end{eg}

\begin{thm}\label{subnot}
Suppose that $P_1\subseteq P$. 
If $P_1$ is not sub-representable, then $P$ is not sub-representable.
\end{thm}

\begin{proof}
Suppose that $P$ was sub-representable and $g$ is as in Definition 
\ref{subrep}. 
Now $g(P_1)$ embeds into $P_1$ by a map $h$. 
Then it is clear that if $g^*=g_{|P_1}$, then $h\circ g^*$ will 
sub-represent $P_1$.
\end{proof}

\begin{defn}
\addtocounter{figure}{1}
The posets whose Hasse diagrams are given in Figure \thefigure \;shall be
known as a {\it vee}, a {\it wedge}, and a {\it diamond}, respectively.
\addtocounter{figure}{-1}
\begin{center}
\begin{picture}(200,100)(0,0) 
\put(30,0){\circle*{10}} 
\put(10,40){\circle*{10}}
\put(50,40){\circle*{10}} 
\put(30,5){\line(-1,2){16}}
\put(30,5){\line(1,2){16}} 
\put(80,0){\circle*{10}} 
\put(100,40){\circle*{10}}
\put(120,0){\circle*{10}} 
\put(100,35){\line(-1,-2){16}}
\put(100,35){\line(1,-2){16}}
\put(170,0){\circle*{10}}
\put(150,40){\circle*{10}} 
\put(190,40){\circle*{10}}
\put(170,80){\circle*{10}}
\put(170,5){\line(1,2){16}}
\put(170,75){\line(-1,-2){16}}
\put(170,75){\line(1,-2){16}}
\put(170,5){\line(-1,2){16}} 
\end{picture}\\
\piccap{}
\end{center}
\end{defn}

\begin{eg}
\mbox{}
\begin{itemize}
\item All two-point posets are (trivially) sub-representable.
\item All three-point posets are (trivially) sub-representable.
\item Not all four-point posets are sub-representable: see Example
\ref{counter}.
\end{itemize}
\end{eg}

\begin{eg}\label{counter}
\addtocounter{figure}{1}
The poset $P$ in Figure \thefigure \;is not sub-representable for the
following reasons. Suppose that $P$ were sub-representable. Since $P$
contains a wedge, the wedge must be sub-represented by either $\{1,3,4\}$
or $\{2,3,4\}$. A two-point chain must then be mapped to $\{1,3\}$,
$\{3,4\}$ or $\{2,3\}$; not to $\{1,2\}$. However, in all three cases we
see that there is no other point incomparable with the chain, and hence
the disjoint union $\{1,2,4\}$ of a two-point chain and a single point
cannot be sub-represented: a contradiction.
\addtocounter{figure}{-1}
\begin{center}
\begin{picture}(100,100)(0,0)
\put(20,10){\circle*{10}}
\put(5,5){$1$}
\put(20,50){\circle*{10}}
\put(5,45){$2$}
\put(20,90){\circle*{10}}
\put(5,85){$3$}
\put(55,50){\circle*{10}}
\put(65,45){$4$}
\put(20,15){\line(0,1){30}}
\put(20,55){\line(0,1){30}}
\put(25,90){\line(2,-3){30}}
\end{picture}\\
\piccap{}
\end{center}

\addtocounter{figure}{1}
Figure \thefigure \;gives all four-point posets that cannot be
sub-represented, as may be verified by simple arguments like that of above.
\addtocounter{figure}{-1}
\begin{center}
\begin{picture}(300,150)(0,0)
\put(10,10){\circle*{5}}
\put(10,12.5){\line(1,1){18}}
\put(30,30){\circle*{5}}
\put(50,10){\circle*{5}}
\put(50,12.5){\line(-1,1){18}}
\put(60,10){\circle*{5}}
\put(130,10){\circle*{5}}
\put(130,12.5){\line(-1,1){18}}
\put(130,12.5){\line(1,1){18}}
\put(110,30){\circle*{5}}
\put(150,30){\circle*{5}}
\put(160,10){\circle*{5}}
\put(210,10){\circle*{5}}
\put(210,12.5){\line(0,1){15}}
\put(230,12.5){\line(0,1){15}}
\put(230,10){\circle*{5}}
\put(210,30){\circle*{5}}
\put(210,12.5){\line(1,1){18}}
\put(230,30){\circle*{5}}
\put(230,12.5){\line(-1,1){18}}
\put(10,80){\circle*{5}}
\put(10,82.5){\line(0,1){15}}
\put(10,100){\circle*{5}}
\put(10,102.5){\line(0,1){15}}
\put(10,120){\circle*{5}}
\put(12.5,120){\line(1,-1){18}}
\put(30,100){\circle*{5}}
\put(70,80){\circle*{5}}
\put(70,82.5){\line(0,1){15}}
\put(70,100){\circle*{5}}
\put(70,102.5){\line(0,1){15}}
\put(70,120){\circle*{5}}
\put(70,82.5){\line(1,1){18}}
\put(90,100){\circle*{5}}
\put(160,80){\circle*{5}}
\put(160,82.5){\line(1,1){18}}
\put(160,82.5){\line(-1,1){18}}
\put(140,100){\circle*{5}}
\put(142.5,100){\line(1,1){18}}
\put(180,100){\circle*{5}}
\put(177.5,100){\line(-1,1){18}}
\put(160,120){\circle*{5}}
\put(220,80){\circle*{5}}
\put(220,82.5){\line(1,1){18}}
\put(240,100){\circle*{5}}
\put(242.5,100){\line(1,-1){20}}
\put(260,80){\circle*{5}}
\put(260,82.5){\line(1,1){18}}
\put(280,100){\circle*{5}}
\end{picture}\\
\piccap{}
\end{center}

\addtocounter{figure}{1}
Figure \thefigure \;gives all four-point posets (likewise
identified) that can be sub-represented.
\addtocounter{figure}{-1}
\begin{center}
\begin{picture}(250,150)(0,0)
\put(15,30){\circle*{5}} \put(25,30){\circle*{5}}
\put(35,30){\circle*{5}} \put(45,30){\circle*{5}}
\put(80,10){\circle*{5}} \put(80,12.5){\line(0,1){15}}
\put(80,30){\circle*{5}} \put(80,32.5){\line(0,1){15}}
\put(80,50){\circle*{5}} \put(80,30){\line(1,1){15}}
\put(95,45){\circle*{5}}

\put(135,10){\circle*{5}} \put(135,12.5){\line(0,1){15}}
\put(135,30){\circle*{5}} \put(135,32.5){\line(0,1){15}}
\put(135,50){\circle*{5}} \put(135,30){\line(1,-1){15}}
\put(150,15){\circle*{5}}

\put(210,30){\circle*{5}} \put(210,30){\line(0,1){15}}
\put(210,30){\line(1,1){15}} \put(210,30){\line(-1,1){15}}
\put(195,45){\circle*{5}} \put(225,45){\circle*{5}}
\put(210,45){\circle*{5}}

\put(30,100){\circle*{5}} \put(30,100){\line(0,-1){15}}
\put(30,85){\circle*{5}} \put(30,100){\line(1,-1){15}}
\put(45,85){\circle*{5}} \put(30,100){\line(-1,-1){15}}
\put(15,85){\circle*{5}}

\put(100,100){\circle*{5}} \put(100,100){\line(0,-1){15}}
\put(100,85){\circle*{5}} \put(110,85){\circle*{5}}
\put(120,85){\circle*{5}}

\put(170,85){\circle*{5}} \put(170,85){\line(0,1){15}}
\put(170,100){\circle*{5}} \put(170,100){\line(0,1){15}}
\put(170,115){\circle*{5}} \put(170,115){\line(0,1){15}}
\put(170,130){\circle*{5}}

\put(220,90){\circle*{5}} \put(220,90){\line(0,1){15}}
\put(220,105){\circle*{5}} \put(220,105){\line(0,1){15}}
\put(220,120){\circle*{5}} \put(230,105){\circle*{5}}
\end{picture}\\
\piccap{}
\end{center}
\end{eg}

Example \ref{counter} shows that we do not have universal 
sub-represent\-ability, even among finite posets. 
We seek to identify which posets are sub-represent\-able.

\begin{note}
Let $x\in P$. We shall use the following notation (note the strictness of
the inequalities):
\begin{itemize}
\item $U(x)=\{y\in P:y>x\}$
\item $D(x)=\{y\in P:y<x\}$.
\end{itemize}
\end{note}

\begin{defn}
A linear ordering $A$ is called a
\begin{itemize}
\item {\it well-ordering} if every non-empty subset $B$ of $A$ has a
least element,
\item {\it well-ordering*} if every non-empty subset $B$ of $A$ has a
greatest element.
\end{itemize}
\end{defn}

\begin{defn}
We shall call $P$ a
\begin{itemize}
\item {\it flower} if and only if there exists $x\in P$ such that $D(x)$
is a well-ordered* chain, $U(x)$ is an antichain (with $|U(x)|>1$) and
$D(x)\cup U(x) \cup \{x\} = P$
\item {\it co-flower} if and only if there exists $x\in P$ such that
$U(x)$ is a well-ordered chain, $D(x)$ is an antichain (with
$|D(x)|>1$) and $D(x)\cup U(x) \cup \{x\}= P$.
\end{itemize}
\end{defn}

\begin{note}
Note that the dual of a flower is a co-flower and vice versa. 
Also every flower contains a vee and every co-flower contains a wedge.
\end{note}

\begin{thm}\label{vandw}
Suppose that $P$ contains both a vee and a wedge. 
Then $P$ is not sub-representable.
\end{thm}

\begin{proof}
Suppose that $P$ is sub-representable but contains both a vee and a wedge. 
Suppose that a vee, $\{1,2,3:1<2,1<3\}$, is represented by
$\{a,b,c:b<a,b<c\}$. 
Then we must have the two-point antichain embedded as $\{a,c\}$. 
The wedge $\{4,5,6:4<6,5<6\}$ would then partially embed as
follows: $\{4,5\} \rightarrow \{a,c\}$, and hence we would have a diamond
$\{a,b,c,d:b<a,b<c,a<d,c<d\}$ contained in $P$. 
By Example \ref{counter}, any space containing a diamond is not 
sub-representable: a contradiction.
\end{proof}

\begin{thm}\label{near-flower}
\addtocounter{figure}{1}
Suppose that $P$ contains one of the posets in Figure \thefigure. 
Then $P$ is not sub-representable.
\addtocounter{figure}{-1}
\begin{center}
\begin{picture}(300,170)(0,0)
\put(45,140){$x$} \put(130,140){$y$}
\put(90,120){\circle*{10}} \put(120,150){\circle*{10}}
\put(90,120){\line(1,1){30}} \put(90,120){\line(-1,1){30}}
\put(60,150){\circle*{10}} 
\put(90,10){\circle*{10}}\put(90,30){\circle*{10}}\put(90,50){\circle*{10}}
\put(90,70){\circle*{10}}
\put(90,10){\line(0,1){20}}\put(90,30){\line(0,1){20}}\put(90,50){\line(0,1){20}}
\put(89,75){.} \put(89,80){.} \put(89,85){.}\put(89,90){.}\put(89,95){.}
\put(10,50){\scriptsize No maximum} 
\put(10,40){\scriptsize element.}
\put(50,30){\vector(1,0){20}}
\put(60,0){\bf A}

\put(190,40){\circle*{10}} \put(220,10){\circle*{10}}
\put(160,10){\line(1,1){30}} \put(190,40){\line(1,-1){30}}
\put(160,10){\circle*{10}} 
\put(190,150){\circle*{10}}\put(190,130){\circle*{10}}\put(190,110){\circle*{10}}
\put(190,90){\circle*{10}}
\put(190,90){\line(0,1){20}}\put(190,110){\line(0,1){20}}\put(190,130){\line(0,1){20}}
\put(189,85){.} \put(189,80){.} \put(189,75){.}\put(189,70){.}\put(189,65){.}
\put(189,60){.} \put(245,0){\bf B}
\put(210,130){\scriptsize No minimum}
\put(210,120){\scriptsize element.}
\put(230,110){\vector(-1,0){20}}
\end{picture}\\
{\rm \piccap{}}
\end{center}
\end{thm}

\begin{proof}
Suppose that $P$ contains poset $A$ but that $P$ is sub-representable. 
Notice that $A$ contains a copy of $\N$. 
Then a two-point antichain must be represented by $\{x,y\}$ and hence a
singleton must be represented by either $\{x\}$ or $\{y\}$. 
In either case we have that the representative of $\N$ would have a
maximal element: a contradiction, and hence $P$ is not
sub-representable. A similar contradiction arises if $P$ contains $B$. 
\end{proof}

\begin{thm}
If $P$ is sub-representable then $P$ is a flower or a co-flower or a
disjoint union of chains.
\end{thm}

\begin{proof}
Suppose that $P$ is neither a flower nor a co-flower nor a disjoint union
of chains. 
Since it is not a disjoint union of chains it contains a vee or a wedge. 
By Theorem \ref{vandw} if it contains both then it is not 
sub-representable. 
Suppose then that $P$ contains a vee but not a wedge.
Since it is not a flower it must contain one of the following:
\begin{center}
\begin{picture}(400,150)(0,0)
\put(50,30){\circle*{10}} \put(80,60){\circle*{10}}
\put(50,30){\line(1,1){30}} \put(50,30){\line(-1,1){30}}
\put(20,60){\circle*{10}} \put(100,40){\circle*{10}} \put(30,0){\bf A}
\put(190,30){\circle*{10}} \put(220,60){\circle*{10}}
\put(220,60){\line(1,1){30}} \put(190,30){\line(1,1){30}}
\put(250,90){\circle*{10}} \put(190,30){\line(-1,1){30}}
\put(160,60){\circle*{10}} \put(190,0){\bf B}
\put(340,120){\circle*{10}} \put(370,150){\circle*{10}}
\put(340,120){\line(1,1){30}} \put(340,120){\line(-1,1){30}}
\put(310,150){\circle*{10}} 
\put(340,10){\circle*{10}}\put(340,30){\circle*{10}}\put(340,50){\circle*{10}}
\put(340,70){\circle*{10}}
\put(340,10){\line(0,1){20}}\put(340,30){\line(0,1){20}}\put(340,50){\line(0,1){20}}
\put(339,75){.} \put(339,80){.} \put(339,85){.}\put(339,90){.}\put(339,95){.}
\put(310,0){\bf C}
\end{picture}\\
\piccap{}
\end{center}

However, $A$ and $B$ are not sub-representable by Example \ref{counter}
and $C$ is not sub-representable by Theorem \ref{near-flower}. 
Hence $P$ is not sub-representable either. Dually, suppose that $P$
contains a wedge but not a vee. Since it is not a
co-flower it must contain one of the following:
\begin{center}
\begin{picture}(400,150)(0,0)
\put(50,60){\circle*{10}} \put(80,30){\circle*{10}}
\put(50,60){\line(1,-1){30}} \put(20,30){\line(1,1){30}}
\put(20,30){\circle*{10}} \put(100,30){\circle*{10}} \put(30,0){\bf D}

\put(190,90){\circle*{10}} \put(220,60){\circle*{10}}
\put(250,30){\line(-1,1){30}} \put(220,60){\line(-1,1){30}}
\put(250,30){\circle*{10}} \put(190,90){\line(-1,-1){30}}
\put(160,60){\circle*{10}} \put(190,0){\bf E}

\put(340,40){\circle*{10}} \put(370,10){\circle*{10}}
\put(310,10){\line(1,1){30}} \put(340,40){\line(1,-1){30}}
\put(310,10){\circle*{10}} 
\put(340,150){\circle*{10}}\put(340,130){\circle*{10}}\put(340,110){\circle*{10}}
\put(340,90){\circle*{10}}
\put(340,90){\line(0,1){20}}\put(340,110){\line(0,1){20}}\put(340,130){\line(0,1){20}}
\put(339,85){.} \put(339,80){.} \put(339,75){.}\put(339,70){.}\put(339,65){.}
\put(339,60){.}\put(339,65){.} \put(340,0){\bf F}
\end{picture}\\
\piccap{}
\end{center}

However, $D$ and $E$ are not sub-representable by Example \ref{counter}
and $F$ is not sub-representable by Theorem \ref{near-flower}. 
Hence neither is $P$.
\end{proof}

\begin{eg}\label{subord}
Let $P = \bigcup_{k\in \N} \{P_k\}$ where $P_k$ is the $k$-point chain. 
Then $P$ is not sub-representable.
\end{eg}

\begin{proof}
Suppose that this poset is sub-representable. 
Then there exist values of $k$ and $m$ such that $P_k$ is represented by
being mapped into $P_m (m\geq k)$. 
We then have that $P_{m+1}$ must map into more than one $P_k$, since it
must contain the image of $P_k$ whereas $P_m$ contains only $m$
points: this yields the desired contradiction.
\end{proof}

\begin{note}
Let $P$ be a {\it finite} poset.
Recall that the {\it height} of $P$, denoted by $ht(P)$, is the largest
cardinality of a chain in $P$, and that the {\it width} of $P$, denoted
by $wd(P)$, is the largest cardinality of an antichain in $P$.
\end{note}

\begin{thm}
Suppose that $P$ is finite. 
If $P$ is either a flower or a co-flower or a disjoint union of chains
then $P$ is sub-representable.
\end{thm}

\begin{proof}
\addtocounter{figure}{1}
Suppose that $P$ is a flower with $ht(P)=n$ and $wd(P)=k$. 
The family of subsets of $P$ consists of flowers of height $m$ and width
$r$ for all $1<m\leq n, 1<r\leq k$, together with chains of size $\leq n$
and antichains of size $\leq k$. 
Label the points of $P$ in the maximal antichain as $x_1,x_2,...,x_k$. 
There exists a point $x$ such that $P= D(x)\cup U(x)\cup \{x\}$. 
Label $x$ as $x_{k+1}$ and the points of $D(x)$ as $x_{k+2},
x_{k+3},...,x_{k+n-1}$ where $x_{k+2}>x_{k+3}>...>x_{k+n-1}$ as in Figure
\thefigure. 
Our isomorphism is defined as follows: map a singleton to $x_1$; map each
chain of size $m>1$ to $\{x_1,x_{k+1},...,x_{k+m-1}\}$; map each
antichain of size $r>1$ to $\{x_1,...,x_r\}$; map each flower of height
$m$ and width $r$ to $\{x_j:k<j\leq k+m-1\}\cup \{x_j:1\leq j\leq r\}$. 
It is clear that $P$ has been sub-represented.
\addtocounter{figure}{-1}
\begin{center}
\begin{picture}(180,160)(0,0)
\put(20,155){$x_1$} \put(80,155){$x_2$} 
\put(170,155){$x_k$}
\put(40,120){\circle*{10}} \put(70,150){\circle*{10}}
\put(160,150){\circle*{10}}
\put(85,150){.} \put(95,150){.}
\put(105,150){.} \put(115,150){.} \put(125,150){.} \put(135,150){.}
\put(145,150){.} \put(155,150){.}
\put(40,120){\line(4,1){120}}
\put(40,120){\line(1,1){30}} \put(40,120){\line(-1,1){30}}
\put(10,150){\circle*{10}} \put(40,80){\circle*{10}} 
\put(40,10){\circle*{10}}
\put(50,10){$x_{k+n-1}$} \put(50,80){$x_{k+2}$}
\put(50,115){$x_{k+1}$} 
\put(40,80){\line(0,1){40}}
\put(39,20){.} \put(39,30){.}
\put(39,40){.} \put(39,50){.} \put(39,60){.} \put(39,60){.}\put(39,70){.} 
\end{picture}\\
\piccap{}
\end{center}

The case where $P$ is a co-flower is precisely dual to the preceding
discussion.
\addtocounter{figure}{1}
Now suppose that $P$ is a disjoint union of chains.
Arrange the chains in descending order of cardinality. 
That is, denote $P$ as $\bigcup \{P_k:k\leq n\}$ where 
$ht(P_k)\geq ht(P_{k+1})$. 
Label the points of each $P_k$ as $\{x^k_m:1\leq m\leq ht(P_k)\}$ as in
Figure \thefigure. 
\addtocounter{figure}{-1}
\begin{center}
\begin{picture}(200,160)(0,0)
\put(20,20){\circle*{10}} \put(20,50){\circle*{10}}
\put(20,80){\circle*{10}} \put(20,110){\circle*{10}}
\put(20,140){\circle*{10}} \put(20,170){\circle*{10}}
\put(20,20){\line(0,1){30}}
\put(20,50){\line(0,1){30}}\put(20,80){\line(0,1){30}}\put(20,140){\line(0,1){30}}
\put(100,20){\circle*{10}} \put(100,50){\circle*{10}}
\put(100,80){\circle*{10}} \put(100,110){\circle*{10}}
\put(100,140){\circle*{10}}
\put(100,20){\line(0,1){30}}
\put(100,50){\line(0,1){30}}
\put(100,110){\line(0,1){30}}
\put(25,20){.} \put(35,20){.}
\put(45,20){.} \put(55,20){.}
\put(65,20){.} \put(75,20){.}
\put(85,20){.} \put(95,20){.}
\put(19,115){.} \put(19,125){.}
\put(19,120){.} \put(19,130){.}
\put(19,135){.} 
\put(99,85){.} \put(99,90){.}
\put(99,95){.} \put(99,100){.}
\put(99,105){.}
\put(95,0){$P_k$}
\put(110,20){$x_1^k$}
\put(110,50){$x_2^k$}
\put(110,80){$x_3^k$}
\put(110,110){$x_{m-1}^k$}
\put(110,140){$x_m^k$}
\put(130,20){.}
\put(140,20){.} \put(150,20){.} 
\put(160,20){\circle*{10}} \put(160,50){\circle*{10}}
\put(160,80){\circle*{10}}
\put(160,20){\line(0,1){30}}
\put(160,50){\line(0,1){30}}
\end{picture}\\
\piccap{}
\end{center}
Let $S$ be a subset of $P$, and suppose that $S$ intersects each of the
posets $\{P_{k_1}, P_{k_2},...,P_{k_t}\}$ where $|S\cap P_{k_q}|\geq
|S\cap P_{k_{q+1}}|$ for all $1\leq q\leq t-1$. For each $r\leq t$ map
$S\cap P_{k_r}$ to $\{x^r_m:m\leq |S\cap P_{k_r}|\}$. 
The poset $P$ has now been sub-represented.
\end{proof}

\begin{cor}
Let $P$ be a finite poset. 
Then $P$ is sub-representable if and only if $P$ is either a flower or a
co-flower or a disjoint union of chains.
\end{cor}

\begin{eg}
Consider $\Z = \{...,-3,-2,-1,0,1,2,3,...\}$, the set of integers with
their usual ordering. Then $\Z$ is not sub-representable.
\end{eg}

\begin{proof}
Suppose that $\Z$ is sub-representable. Let the image of the equivalence
class of the natural numbers, i.e. $\theta([\N])$, be the set
$\{0',1',2',3',...\}$. Let $\theta([-\N])$ be the set
$\{(-0)',(-1)',(-2)',(-3)',...\}$. Since a singleton must be mapped into
both $\theta([\N])$ and $\theta([-\N])$, we know that $|\theta([\N])\cap
\theta([-\N])|>0$. Suppose that $p'$ is the least point of
$\theta([\N])\cap\theta([-\N])$, that is, $p'=(-q)',$ which implies that
$|\theta([\N])\cap \theta([-\N])|\leq q+1$. Then
$\theta([\{0,1,2,3...,q,q+1\}])$ must be a subset of $\theta([\N])$ and
also a subset of $\theta([-\N])$. This is a contradiction, since
$|\theta([\N])\cap \theta([-\N])|\leq q+1$.
\end{proof}

\begin{cor}
By Corollary \ref{subnot} we then have that $\Q, \R$ and $\R\setminus \Q$
(as ordered sets) cannot be sub-represented.
\end{cor}

\begin{cor}
If $P$ contains a copy of $\N$ and a copy of $-\N$ then it cannot be
sub-represented.
\end{cor}

\begin{note}
If $P$ is a chain which is not well-ordered then it contains a copy of
$-\N$, and if $P$ is a chain which is not well-ordered* then it contains
a copy of $\N$.
\end{note}

\begin{cor}
If $P$ is a chain which is not well-ordered and not well-ordered* then it
cannot be sub-represented.
\end{cor}

\begin{thm}
Every well-ordered set can be sub-represented.
\end{thm}

\begin{proof}
Let $P$ be a well-ordered chain of order-type $\alpha$. It suffices to
show that $\alpha$ itself is sub-representable. Let $S\subseteq \alpha$
and suppose that $S$ is order-isomorphic to $\beta<\alpha$. We map $[S]$
to $\beta$ where $\beta$ is a proper initial segment of $\alpha$. On the
other hand, if $S$ is order-isomorphic to $\alpha$, map $[S]$ to $\alpha$
itself.
\end{proof}

\begin{thm}\label{repwell*}
Every well-ordered* set can be sub-represented.
\end{thm}

\begin{proof}
Let $P$ be a well-ordered* chain. Let $P^*$ denote the dual of $P$. Then
$P^*$ is, without loss of generality, an ordinal $\alpha$. If $S\subseteq
P$ then $S^*$ is order-isomorphic to some $\beta\leq\alpha$, and so map
$[S]$ to the final segment of $\alpha^*$ of order-type $\beta^*$.
\end{proof}

\begin{thm}
Let $P$ be a flower or a co-flower. Then $P$ is sub-represent\-able.
\end{thm}

\begin{proof}
Suppose that $P$ is a flower. Then $P= D(x)\cup U(x) \cup \{x\}$. Suppose
that $|U(x)|=\delta$. Label the points of $U(x)$ as
$\{x_\beta:\beta<\delta\}$. By definition, $D(x)\cup \{x\} \cup \{x_0\}$
is a well-ordered* chain. Label the points of $U(x)\setminus \{x_0\}$ as
$\{x'_\beta:\beta<\delta'\}$. Let $S\subseteq P$. If $S$ is a copy of
$\alpha^*$ for some ordinal $\alpha$, map $S$ to the final segment
$\alpha^*$ of $D(x)\cup \{x\} \cup \{x_0\}$. Otherwise $S$ contains a copy
of $\alpha^*$ for some ordinal $\alpha$ (such that it does not contain a
copy of $\gamma^*$ for some $\gamma>\alpha$) together with an antichain of
cardinality $\zeta$. Map $S$ to $\{x'_\beta:\beta<\zeta\}$ together with
the final segment $\alpha^*$ of $D(x)\cup \{x\} \cup \{x_0\}$. We have now
sub-represented $P$. The case where $P$ is a co-flower is similar.
\end{proof}

\begin{thm}\label{finord}
Let $P$ be the disjoint union of finitely many well-ordered sets.
Then $P$ is sub-representable.
\end{thm}

\begin{proof}
Suppose that $P= \bigcup_{k\leq n} \bigcup_{s\leq m_k} \{^s\omega^k\}$
where $^s\omega^k$ is a copy of an ordinal $\omega^k$. Then (up to
order-isomorphism) $P$ is a copy of $n^*=\sum_k m_k$ many ordinals. So let
us consider $P$ as $\bigcup_{n\leq n^*}\{\omega^n\}$ where $\omega^n\geq
\omega^{n+1}$ for all $n\leq n^*$. Let $Y\subseteq P$ meet each of
$\omega^{n_1},\omega^{n_2},...,\omega^{n_t}$ where $|Y\cap
\omega^{n_p}|\geq |Y\cap \omega^{n_{p+1}}|$ for all $1\leq p\leq t-1$. Now
$Y\cap \omega^{n_p}$ is a copy of an ordinal $\alpha_p$. So map $Y\cap
\omega^{n_p}$ to $\alpha_p\in \omega^p$ and hence $P$ has been
sub-represented.
\end{proof}

\begin{cor}
A similar proof gives the result that a disjoint union of finitely many
well-ordered* sets is sub-representable.
\end{cor}

Let us examine further the question of which disjoint unions of ordinals
can be sub-represented. Let $P$ be a disjoint union of ordinals. If $P$
contains a copy of the disjoint union of the family of ordinals
$\{1,2,3,4,...\}$ then by Example \ref{subord} and Corollary
\ref{subnot}, $P$ cannot be sub-represented. Hence $P$ must contain
copies of only finitely many distinct ordinals and only finitely many
copies of each relevant infinite ordinal. If we have only finitely many
copies of each ordinal then by Theorem \ref{finord}, $P$ is
sub-representable. So all that remains to consider is the case where
there are infinitely many copies of one or more finite ordinals. It will
be seen --- although the demonstration requires considerably more effort
than earlier proofs in this paper --- that this case also is
sub-representable.

\section{Sub-representation of pinboards}

\begin{defn}
A {\it pinboard} is a finite set of ordered pairs 
$$\{(h_i,f_i): 1\leq i\leq k\}$$ 
in which, for each value of $i$, $h_i$ (the height) is an ordinal, $f_i$
(the frequency) is a cardinal, and not both $h_i,f_i$ can be infinite. 
We recall that a cardinal is merely an initial ordinal (or equivalently,
the least one of a given cardinality).
\end{defn}

\begin{defn}
A {\it simple pinboard} is a pinboard of the form
$$\{(\beta,n),(m,\gamma)\}$$ where $\beta$ and $\gamma$ are infinite
cardinals, and $m$ and $n$ are finite cardinals.
\end{defn}

\begin{defn}
\addtocounter{figure}{1}
The {\it poset of a pinboard} $\{(h_i,f_i): 1\leq i\leq k\}$ is the
disjoint union of $f_i$-many copies of $h_i$ for $1\leq i\leq k$. 
For example, if the pinboard was 
$\{(\omega_2,5),(\omega_1,2),(6,\omega),(3,1)\}$, then the poset is as
suggested by Figure \thefigure.
\addtocounter{figure}{-1}
\end{defn}

\begin{center}
\begin{picture}(300,200)(0,0)
\put(7,175){$\omega_2$}
\put(7,115){$\omega_1$}
\put(10,75){$6$}
\put(10,26){$3$}
\put(20,20){\line(1,0){240}}
\put(20,20){\line(0,1){160}}
\put(40,20){\line(0,1){160}}
\put(60,20){\line(0,1){160}}
\put(80,20){\line(0,1){160}}
\put(100,20){\line(0,1){160}}
\put(120,20){\line(0,1){100}}
\put(140,20){\line(0,1){100}}
\put(160,20){\line(0,1){60}}
\put(180,20){\line(0,1){60}}
\put(260,20){\line(0,1){10}}
\put(18,10){$0$}
\put(38,10){$1$}\put(58,10){$2$}\put(78,10){$3$}
\put(98,10){$4$}\put(118,10){$5$}\put(138,10){$6$}
\put(158,10){$7$}\put(178,10){$8$}
\put(258,10){$\omega_0$}
\put(190,10){.}\put(200,10){.}\put(210,10){.}\put(220,10){.}
\put(230,10){.} \put(240,10){.}\put(250,10){.}
\end{picture}\\
\piccap{}
\end{center}

Our aim is to sub-represent the poset of an arbitrary pinboard.
However, it is easy to see that the poset of any pinboard is a subset
of the poset of a simple pinboard. Hence, by Theorem \ref{subnot},
if every simple pinboard's poset is sub-representable, then every
pinboard's poset is sub-representable.

Let $P$ be the poset of a simple pinboard $\{(\beta,n),(m,\gamma)\}$. 
We shall represent $P$ by the set
$$X=\{(\beta,i,\alpha):\alpha<\beta, i<n\}\cup \{(m,\alpha,j):j<m,
\alpha<\gamma\},$$ 
that is, an isomorphic copy of $\beta$ over each $(\beta,i)$ and an
isomorphic copy of $m$ over each $(m,\alpha)$. We note that the ordinal
$\oplus_m \gamma$ has the same cardinality as the set
$$F=\{(\beta,i):i<n\}\cup \{(m,\alpha):\alpha<\gamma\}$$ 
which underlies $X$. 
We can therefore choose an injection $\lambda:\oplus_m\gamma\rightarrow F$
such that the first $n$ terms of its domain map onto 
$\{(\beta,i):i<n\}$. 
We now have labelled the columns of $X$ in a $\oplus_m \gamma$-sequence
such that the first $n$ columns are those of the infinite ordinals.

Let $Y$ be a subset of the simple pinboard above. 
Then $Y$ consists of at most $n$ ordinals that exceed $m$ and at most 
$\gamma$-many copies of each of the ordinals $m,m-1,...,2,1$. 
We can assume that for each term $(h,f)$ in $Y$ there does not exist a
term $(h',f')$ in $Y$, where $f'$ is infinite, such that $h'<h$ and
$f'<f$: for, if there was such a term $(h',f')$, the elimination of
$(h',f')$ would not affect the embeddability class of $Y$. 
Therefore we can consider $Y$ to be of the following form:
$$Y=
\{(\beta_1,n_1),(\beta_2,n_2),...,(\beta_k,n_k),(m,\gamma_m),
(m-1,\gamma_{m-1}),...,(1,\gamma_1)\},$$ 
where $\beta_1>\beta_2>...>m$.

We note that the ordinal sum 
$$n_1 \oplus n_2 \oplus \cdots \oplus n_k \oplus \gamma_m \oplus \cdots 
\oplus \gamma_1 \leq 
n \oplus(\oplus_m \gamma) = 
\oplus_m \gamma$$ 
and hence we can find a `remainder' ordinal $\gamma_0$ such that
$$n_1 \oplus n_2 \oplus \cdots \oplus n_k \oplus \gamma_m \oplus \cdots
\oplus \gamma_1 \oplus \gamma_0 = \oplus_m \gamma.$$ 
We now apply our injection $\lambda$ and hence associate a copy of 
$\beta_1$ with the first $n_1$-many columns of $X$, a copy of $\beta_2$
with the next $n_2$-many columns and so on.
Our subset $\theta(Y)$ of $X$ associated with $Y$ is now formed by taking
the desired initial segments of the associated columns, and it is clear
that $\theta(Y)$ is an isomorphic copy of $Y$. 

\begin{thm}
Let $Y,Y'$ be subsets of $X$. 
Then $Y\hookrightarrow Y'$ if and only if $\theta(Y)\subseteq \theta(Y')$.
\end{thm}

\begin{proof}
Suppose that $Y\hookrightarrow Y'$. 
Let us denote $Y'$ as
$$\{\beta_1',n_1'),(\beta_2',n_2'),\ldots,(\beta_l',n_l'),(m,\gamma_m'),
(m-1,\gamma_{m-1}'),\ldots,(1,\gamma_1')\},$$ 
where $\beta_1'>\beta_2'>\cdots>\beta_l'>m$. 
Let $\theta(Y')$ be generated as above. 
We proceed to show that $\theta(Y)\subseteq \theta(Y')$. 
Consider the ordinals (without repetitions) comprising $Y$ and $Y'$, and
denote this set arranged in decreasing order as
$$\{\zeta_1,\zeta_2,\ldots,\zeta_q,m,m-1,\ldots,1\}.$$
It is clear that a single column of $Y$ must embed into a single column of
$Y'$ and hence the ordinal $\zeta_1$ is in fact $\beta_1'$. 
We can now assume that $Y$ and $Y'$ are of the following form where some
of the frequencies may be zero:
$$Y = \{\zeta_1,r_1),(\zeta_2,r_2),\ldots,(\zeta_q,r_q),(m,\gamma_m),
(m-1,\gamma_{m-1}),\ldots,(1,\gamma_1)\},$$
$$Y' = \{\zeta_1,r_1'),(\zeta_2,r_2'),\ldots,(\zeta_q,r_q'),(m,\gamma_m'),
(m-1,\gamma_{m-1}'),\ldots,(1,\gamma_1')\}.$$
The proof now proceeds as follows:

\noindent
Stage 1: 
Any column of $Y$ of height $\zeta_1$ must embed into a column of $Y'$ of
that height, and so it follows that $r_1\leq r_1'$. 
That is to say, the columns of $\theta(Y)$ of height $\zeta_1$ are
contained in $\theta(Y')$.

\noindent
Stage 2: 
Any column of height $\zeta_1$ or $\zeta_2$ must embed into a column of
$Y'$ of height $\zeta_1$ or $\zeta_2$, and hence 
$r_1+r_2\leq r_1'+r_2'$. 
It follows that the columns of $Y$ of height $\zeta_2$ are contained in
$\theta(Y')$.

\noindent
The obvious iteration process will now yield, in a finite number of
stages, $\theta(Y)\subseteq\theta(Y')$.
\end{proof}

\begin{eg}
Suppose that $P$ is the poset of the simple pinboard
{$\{$}{$(\aleph_2,12)$}, {$(7,\aleph_3)$}{$\}$} and we have a labelling
map $\lambda:\oplus_7 \aleph_3 \rightarrow F$ where
$F=\{(\aleph_2,i):i<12\}\cup\{(7,\alpha):\alpha<\aleph_3\}$, 
such that $\{\lambda(0),..,\lambda(11)\}$ are the base points of the 
infinite columns. 
Let $Y$ be a subset of $P$ whose columns have the following ordinalities:
$\{\omega_1\oplus1, \omega_1, \omega_0\oplus5, \omega_0\oplus 5,
\omega_0, 30, 30, 20\}$ 
together with $\aleph_0$ copies of $5$ and $\aleph_0$ copies of $3$. 
As before, we can disregard the $\aleph_0$ copies of $3$ since these do
not affect the embeddability class of $Y$. 
The remaining frequencies (in decreasing order of corresponding
height) are as follows: $\{1,1,2,1,2,1,\aleph_0\}$. 
We then have that $\theta(Y)$ consists of the following initial
segments:
$$\theta(Y) =
\begin{cases}
\omega_1\oplus1 & {\rm of \;column\;} \lambda(0),\\
\omega_1 & {\rm of \;column\;} \lambda(1),\\
\omega_0\oplus5 & {\rm of \;columns\;} \lambda(2) {\rm \;and\;} \lambda(3),\\
\omega_0 & {\rm of \;column\;} \lambda(4),\\
30 & {\rm of \;columns\;} \lambda(5) {\rm \;and\;} \lambda(6),\\
20 & {\rm of \;column\;} \lambda(7),\\
5 & {\rm of \;columns\;} \lambda(8\oplus t) {\rm \;for\; all\;} t<\aleph_0,\\
0 & {\rm otherwise.}
\end{cases}$$

Let $Y'$ be another subset of $P$ whose columns have the following
ordinalities:
$\{\omega_2, \omega_2, {\omega_1\oplus10}, \omega_1, \omega_0, 60, 40, 30,
20\}$ 
together with $\aleph_1$ copies of $6$. 
We note that $Y\hookrightarrow Y'$. 
The frequencies (in decreasing order of corresponding height) are as
follows: $\{2,1,1,1,1,1,1,1,\aleph_1\}$. 
We then have that $\theta(Y')$ consists of the following initial
segments:
$$\theta(Y') =
\begin{cases}
\omega_2 & {\rm of \;columns\;} \lambda(0) {\rm \;and\;} \lambda(1),\\
\omega_1\oplus10 & {\rm of \;column\;} \lambda(2),\\
\omega_1 & {\rm of \;column\;} \lambda(3),\\
\omega_0 & {\rm of \;column\;} \lambda(4),\\
60 & {\rm of \;column\;} \lambda(5),\\
40 & {\rm of \;column\;} \lambda(6),\\
30 & {\rm of \;column\;} \lambda(7),\\
20 & {\rm of \;column\;} \lambda(8),\\
6 & {\rm of \;columns\;} \lambda(9\oplus t) {\rm \;for\; all\;} t<\aleph_1,\\
0 & {\rm otherwise,}
\end{cases}$$
and we observe that $\theta(Y)\subseteq \theta(Y')$.
\end{eg}

\begin{defn}
A {\it co-pinboard} is a finite set of ordered pairs
$\{(h_i,f_i): 1\leq i\leq k\}$ such that $\{(h_i^*,f_i): 1\leq i\leq k\}$ 
is a pinboard.
\end{defn}

\begin{res}
Let $P$ be a poset. 
Then $P$ is sub-representable if and only if $P$ is either:
\begin{itemize}
\item[(i)] a flower, or
\item[(ii)] a co-flower, or
\item[(iii)] the poset of a pinboard, or
\item[(iv)] the poset of a co-pinboard.
\end{itemize}
\end{res}

\begin{note}
We could similarly define sub-representation of an arbitrary topological
space as follows: 
Let $(X,\mathcal{T})$ be a topological space. 
We shall say that $(X, \mathcal{T})$ is {\it sub-representable} if
and only if there exists a map $g$ from ${\mathbb P}(X)$ to 
${\mathbb P}(X)$ such that for all $Y$, $Z \in {\mathbb P}(X)$:
\begin{itemize}
\item[(i)] 
$(Y,\mathcal{T}_{|Y}) \hookrightarrow (Z,\mathcal{T}_{|Z})$ if and
only if 
$g(Y)\subseteq g(Z)$, and
\item[(ii)] 
$(Y,\mathcal{T}_{|Y})\hookrightarrow
(g(Y),\mathcal{T}_{|{g(Y)}}) \hookrightarrow
(Y,\mathcal{T}_{|Y})$.
\end{itemize}
\end{note}

It then follows that we have already characterised all principal $T_0$
spaces which are sub-representable. The question of which other
topological spaces are sub-representable remains open.

\providecommand{\bysame}{\leavevmode\hbox to3em{\hrulefill}\thinspace}
\providecommand{\MR}{\relax\ifhmode\unskip\space\fi MR }
\providecommand{\MRhref}[2]{%
  \href{http://www.ams.org/mathscinet-getitem?mr=#1}{#2}
}
\providecommand{\href}[2]{#2}

\end{document}